\documentclass[10pt, a4paper]{article}
\usepackage[utf8]{inputenc}
\usepackage[english]{babel}
\usepackage{graphicx}
\usepackage{wrapfig}
\usepackage{lmodern}
\usepackage{multicol}
\usepackage{blindtext}
\usepackage{lipsum}
\usepackage{amsmath}
\usepackage{amssymb}

\usepackage{amsfonts}
\usepackage{enumitem}
\usepackage{geometry}
\usepackage[symbol]{footmisc}

\usepackage[font=footnotesize,labelfont=bf]{caption}
\usepackage{subcaption}
\usepackage[figurename=Fig.]{caption}
\usepackage{multicol}
\usepackage{setspace}
\usepackage{indentfirst}
\usepackage{cite}
\usepackage{float}
\usepackage{titlesec}
\titlelabel{\thetitle.}
\usepackage{hyperref}
\hypersetup{
    colorlinks=true,
    citecolor=red,
    linkcolor=blue,
    urlcolor=cyan,,
    }
\newtheorem{thm}{Theorem}[section]

\newtheorem{lem}[thm]{Lemma}
\newtheorem{prop}[thm]{Proposition}
\newtheorem{defn}[thm]{Definition}
\newtheorem{rem}[thm]{\bf{Remark}}
\newtheorem{exm}[thm]{Example}

\numberwithin{equation}{section}

\begin{document}
\begin{center}
\textbf{\Large{Solvability of Infinite Systems of Nonlinear Caputo Fractional Differential equations in  generalized Hahn Sequence Space}}
\end{center}
\begin{center}
Khurshida Parvin$^{a,}$\footnote{e-mail: khurshidaparvin473@gmail.com}, Bipan Hazarika$^{a,}$\footnote{e-mail: bh\_rgu@yahoo.co.in (Corresponding author)} and Awad A. Bakery$^{b,c}$ \footnote{e-mail: awad\_bakery@yahoo.com}\\
$^{a}$Department of Mathematics, Gauhati University, Guwahati-781014, Assam, India\\
$^b$University of Jeddah, College of Science and Arts at Khulis, Department of Mathematics, Jeddah, Saudi Arabia\\
$^c$Ain Shams University, Faculty of Science, Department of Mathematics, Cairo, Abbassia, Egypt 
\end{center}
\textbf{Abstract:}
This paper presents the Hausdorff measure of noncompactness (MNC) within the framework of the generalized Hahn sequence space. By applying the MNC, we explore the existence of solutions for nonlinear Caputo fractional differential equations subject to three-point integral boundary conditions in the generalized Hahn sequence space. We apply certain sufficient conditions to ensure the uniqueness of solutions to the aforementioned problem, utilizing the Banach fixed-point theorem, and discuss the Hyers-Ulam stability of the problem. In conclusion, the analytical framework is complemented by illustrative examples that demonstrate the validity and applicability of the main results.\\



\section{ Introduction}
In recent years, numerous researchers have explored boundary value problems associated with nonlinear fractional differential equations. Fractional derivatives work as a powerful method for describing the memory and hereditary characteristics of various materials and processes \cite{I}. In recent years, numerous researchers have explored boundary value problems associated with nonlinear fractional differential equations. Fractional derivatives serve as a powerful tool for describing the memory and hereditary characteristics involves  diverse materials and methods. The unique properties of fractional derivatives render fractional-order models more applicable and accurate than traditional integer-order models. Within diverse fields of scientific and engineering e.g., physics, chemistry, economics, blood flow studies and experimental data modeling, fractional differential equations frequently appear \cite{I,AA}. Employing the Caputo fractional derivative ensures that the model remains mathematically sound while being practically relevant to the phenomena being analyzed. 

The measure of noncompactness (MNC) serves as a fundamental tool in the analysis of nonlinear functional problems within Banach spaces. In metric fixed-point theory and the study of operator equations in Banach spaces, they play a vital role, particularly in characterizing different types of compact operator. In examining different forms of differential and integral equations, MNC are also applied. Broadly defined, the function measure of noncompactness is defined on all nonempty bounded subsets of a given metric space, which equals zero only for subsets that are relatively compact.

The concept of measure noncompactness represented by $\alpha$, was initially introduced and analyzed by Kuratowski in 1930 \cite{Kura}. Later in 1955, the symbol $\alpha$ was first applied by Darbo \cite{Darbo} in proving the fixed-point theorem, which is useful in investigating existence results for different classes of operator equations.

The second measure of noncompactness, known as  Hausdorff or ball measure, is denoted by $\chi$ and it was initiated by Golden\v{s}tein et al. in 1957 \cite{G} and subsequently examined by Golden\v{s}tein and Markus \cite{ma}), the inner Hausdorff of noncompactness as well as  Istr\v{a}tescu  \cite{I1988}.

In recent studies, the use of MNC in sequence spaces has been investigated for different kinds of differential equations \cite{Alotaibi,Bipan,A.Das,S.Deb,H.Jafari,M.Mursaleen,MM,Mm,M.Mur}.

Farahi et al. \cite{Farahi} examined the existence of solutions for infinite systems of fractional equations in the Hahn sequence space.

Aghajani et al. \cite{Agha}, studied the existence of solutions for infinite systems of second-order differential equations in $\ell_{1}$ spaces.  

To the best of our knowledge, no research has been conducted on fractional differential equations in the generalized Hahn space. The generalized Hahn sequence space expands on classical sequence spaces, creating a more inclusive framework for analysis. Therefore, we are motivated to investigate the existence, uniqueness and HU-stability of infinite systems for nonlinear fractional differential equations (NFDE) with order $\beta\in(0, 1]$
\begin{equation}\label{1.1}
    \begin{split}
        ^c\mathfrak{D}^{\beta}\mathfrak{m}(\xi)=\varphi_{i}(\xi,~\mathfrak{m}(\xi)),~0<\xi<1,~0<\beta\leq1,~i\in\mathbb{N},\\
        \mbox{with~B.C.~}
        \mathfrak{m}(0)=0,~\mathfrak{m}(1)=\mu\int\limits_{0}^{\varrho}\mathfrak{m}(s)ds,~ 0<\varrho<1
    \end{split}
\end{equation}
in the generalized Hahn sequence space, where $^c\mathfrak{D}^{\beta}$ represents the Caputo fractional derivative with order $\beta$. The function $\varphi_{i}:[0,1]\times\mathbb{R} \rightarrow \mathbb{R}$ is continuous and $\mu\neq\frac{2}{\varrho^{2}}$.

Let $(\Omega, ~\lVert.\rVert)$ be a real Banach space with zero element 0 and denote by $\mathfrak{D}(\mathfrak{m}, \mathfrak{r})$, the closed ball centered at $\mathfrak{m}$ with radius $\mathfrak{r}.$ For any nonempty subset 
$\mathcal{R}\subseteq \Omega$, let $\Bar{\mathcal{R}}$ represents its closure and $Conv(\mathcal{R})$ denote its closed convex hull. The space $\Omega$ contains two notable collections: one comprising all nonempty bounded subsets, denoted by $\mathcal{N}_{\Omega}$ and another consisting of all subsets that are relatively compact, referred to as $\mathcal{M}_{\Omega}$.

\begin{defn}\cite{A}
    A function $\psi: \mathcal{N}_{\Omega} \rightarrow [0, \infty]$ is called a measure of non-compactness (MNC), if it satisfies the following conditions
    \begin{enumerate}
    \item [(1)] The set $ker\psi =\big\{\mathcal{R}\in \mathcal{N}_{\Omega};~ \psi(\mathcal{R})=0\big\}$ is nonempty and $ker\psi \subseteq \mathcal{M}_{\Omega}.$
    \item[(2)] $\mathcal{R}\subseteq Z\implies \psi(\mathcal{R})\leq \psi(Z)$.
    \item[(3)] $\psi(\overline{\mathcal{R}})=\psi(\mathcal{R})$.
    \item[(4)] $\psi(Conv(\mathcal{R}))=\psi(\mathcal{R})$.
    \item[(5)] $\psi(\theta \mathcal{R} + (1-\theta)Z)\leq \theta\psi(\mathcal{R})+(1-\theta)\psi (Z)$ for $\theta \in [0, 1]$.
    \item[(6)] If ${\mathcal{R}_{j}}$ represents a sequence of closed sets from $\mathcal{N}_{\Omega}$ such that for $\mathcal{R}_{j+1}\subseteq \mathcal{R}_{j}$ and $\lim\limits_{j\to\infty}\psi(\mathcal{R}_{j})=0$ for $j=1,2,...$. So, $\bigcap_{j=1}^{\infty}\mathcal{R}_{j}\neq \emptyset$.
    \end{enumerate}
\end{defn}
\begin{defn}\cite{Banas}
    Consider a metric space $(\mathcal{R},d)$ and $\mathcal{A}\in \Omega_{\mathcal{R}}$. The Kuratowski measure of non-compactness of $\mathcal{A}$ is introduced as follows:
    $$\alpha(\mathcal{A})=\inf\big(\varepsilon>0: \bigcup\limits_{i=1}^{p}S_{i}\supseteq\mathcal{A}, ~S_{i}\subset\mathcal{R},~ \varepsilon<diam(S_{i}),~i\in\mathbb{N}\big);$$
    where $diam(S_{i})=\sup\{d(\mathfrak{y},\mathfrak{m}); \mathfrak{y},\mathfrak{m}\in S_{i}\}$ ~and $(i=1,2,....p;~ p\in\mathbb{N}).$\\
    
    The Hausdorff measure of noncompactness $\chi(\mathcal{A})$ on closed ball $\mathfrak{D}(\mathfrak{m}_{i}, \mathfrak{r}_{i})$:
    $$\chi(\mathcal{A})=\inf\big(0<\varepsilon: \mathcal{A}\subset\bigcup\limits_{i=1}^{p}\mathfrak{D}(\mathfrak{m}_{i}, \mathfrak{r}_{i}),~\mathfrak{m}_{i}\in\mathcal{R}, ~\mathfrak{r}_{i}<\varepsilon,~i=1,2,...p\in\mathbb{N}\big).$$
\end{defn}
\begin{defn}\cite{M}
   Suppose $\mathcal{Q}$ be a subset of a Banach space $\Omega$ where $\mathcal{Q}\neq\phi$ and let $\psi$ be a MNC on $\Omega$. Also $\mathcal{R}$ is a bounded subset of $\mathcal{Q}$. For a mapping $\mathcal{H}:\mathcal{Q}\rightarrow \mathcal{Q}$ is a Meir-Keeler condensing (MKC) operator
   if for every $0<\varepsilon$, there exists $0<\delta$ such that 
   for any subsets $\mathcal{R}\subseteq\mathcal{Q},$ the condition 
   $$\varepsilon\leq \psi(\mathcal{R})<\delta+\varepsilon ~implies~ \psi(\mathcal{H}(\mathcal{R}))<\varepsilon$$ is satisfied.
\end{defn}
\begin{thm}\cite{M}\label{1.4}
    Consider a Banach space  $\Omega$ and let $\mathfrak{D}\subseteq \Omega$ be a convex set that is bounded and nonempty. Assume $\psi$ is a measure of non-compactness defined on $\Omega$. If the mapping $\mathcal{H}:\mathfrak{D}\rightarrow\mathfrak{D}$ satisfies the Meir-Keeler condensing operator, then $\mathcal{H}$ has a fixed point.
\end{thm}
\begin{prop}\cite{Banas}\label{1.5}
    Let $\mathfrak{Y}$ be a bounded and equicontinuous subset of $\mathcal{C}(\kappa, \Omega)$ and let $\chi$ denote a measure of noncompactness. Then, the mapping  $\chi(\mathfrak{Y}(.))$ is continuous and the following inequalities hold:
    $$\sup\limits_{\xi\in\kappa}\chi(\mathfrak{Y}(\xi))=\chi(\mathfrak{Y}).\\
    \mbox{~Also~}~\chi\left(\int\limits_{0}^{\xi}\mathfrak{Y}(\kappa)d\kappa\right)\leq\int\limits_{0}^{\xi}\chi(\mathfrak{Y}(\kappa))d\kappa.$$
\end{prop}

\begin{thm}\cite{Smart}(Banach Fixed Point Theorem)\label{1.7}:
    If $\mathcal{H}$ is a contraction mapping in a Banach space $\Omega$. Then $\mathcal{H}$ has unique fixed point.
\end{thm}

\begin{defn}\cite{AA}
     The Caputo derivative of fractional order $\beta$ for the continuous function $\varphi:[0,\infty)\rightarrow \mathbb{R}$,  where $n-1<\beta<n$ and $n=[\beta]+1,$ is expressed as 
    $$^\mathcal{C}\mathfrak{D}^{\beta}\varphi(\xi)=\frac{1}{\Gamma(n-\beta)}\int\limits_{0}^{\xi}(\xi-s)^{n-\beta-1}\varphi^{(n)}(s)ds,$$
    where $[\beta]$ represents the integer part of the real number $\beta$.
\end{defn}    
\begin{defn}\cite{SG}
   The Riemann-Liouville fractional integral of order $\beta$ at point $\xi$ is given by 
   $$I^{\beta}\varphi(\xi)=\frac{1}{\Gamma(\beta)}\int\limits_{0}^{\xi}\frac{\varphi(s)} {(\xi-s)^{1-\beta}}ds,~~\beta>0,$$
   for a function $\varphi$ defined on $[0,\infty)$ and a real number and assuming the integral is well-defined.
\end{defn}
\begin{defn}\cite{SG}
    The Riemann-Liouville fractional derivative of order $\beta$ is as:
    $$\mathfrak{D}^{\beta}\varphi(\xi)=\frac{1}{\Gamma(n-\beta)}\bigg(\frac{d}{d\xi}\bigg)^{n}\int\limits_{0}^{\xi}\frac{\varphi(s)}{(\xi-s)^{\beta-n+1}}ds,$$
    provided that the integral on the right-hand side exists for all $\xi>0$ and  $\varphi$ is a continuous function on the interval $(0,\infty)$ and $\beta>0$ with $n=[\beta]+1$.
\end{defn}

\begin{lem}\cite{Bashir}
    Consider the boundary value problem defined in the Equation (\ref{1.1}). Its solution $\mathfrak{m}(\xi)$ can be represented as
    \begin{equation*}
        \begin{split}
           \mathfrak{m}(\xi)
           &=\frac{1}{\Gamma(\beta)}\int\limits_{0}^{\xi}(\xi-s)^{\beta-1}\varphi(s,\mathfrak{m}(s))ds-\frac{2\xi}{(2-\mu\varrho^{2})\Gamma({\beta})}\int\limits_{0}^{1}(1-s)^{\beta-1}\varphi(s,\mathfrak{m}(s))ds\\
           &+\frac{2\mu\xi}{(2-\mu\varrho^{2})\Gamma({\beta})}\int\limits_{0}^{\varrho}\left( \int\limits_{0}^{s}(s-\mathfrak{n})^{\beta-1}\varphi(\mathfrak{n},\mathfrak{m}(\mathfrak{n})d\mathfrak{n}\right)ds.
        \end{split}
    \end{equation*}
\end{lem}

\section{ Generalized Hahn Sequence Space}
Let $\omega$ represent the set of all complex sequences $\mathfrak{m}=(\mathfrak{m}_{n})_{n=1}^{\infty}$. Additionally, 
the sets of bounded sequences, convergent sequences, null sequences and finite sequences (those terminating in zeroes) are represented by $\ell_{\infty}$, $c$,~ $c_{0}$ and $\phi$ respectively. Moreover, $cs$, $bs$ and $\ell^{1}$ represent the sets of all convergent series, bounded series and absolutely convergent series, respectively. We denote $e=(e_{n})_{n=1}^{\infty}$ and $e^{(j)}=(e_{n}^{(j)})_{n=1}^{\infty}$ $(j\in\mathbb{N})$ for the sequences with $e_{n}=1$ for all $n$, and $e_{j}^{(j)}=1$ and $e_{n}^{(j)}=0$ for $n\neq j$.

An FK-space is a type of complete linear metric sequence space in which every convergent sequence also converges point-wise in each coordinate. Also, if an FK-space is equipped with a norm, it is called a BK-space. A BK-space $\mathfrak{X}$ is called an AK-property if, for any sequence $\mathfrak{m}=(\mathfrak{m}_{n})_{n=1}^{\infty}\in\mathfrak{X}$, we have  $\mathfrak{m}=\lim\limits_{r\rightarrow\infty}\mathfrak{m}^{[r]}$, where the $r$-section of the given sequence $\mathfrak{m}$ is represented as $\mathfrak{m}^{[r]}=\sum\limits_{n=1}^{r}\mathfrak{m}_{n}e^{(n)}$. Some BK-spaces are $\ell_{\infty}$, $c$, $c_{0}$, $cs$, $bs$ and $\ell_{1}$ according  to their norms $\lVert\mathfrak{m}\rVert_{bs}=\sup\bigg|\sum\limits_{n=1}^{j}\mathfrak{m}_{n}\bigg|$ on $cs$ and $bs$, $\lVert\mathfrak{m}\rVert_{1}=\sum\limits_{n=1}^{\infty}|\mathfrak{m}|$ on $\ell_{1},$ $\Vert\mathfrak{m}\rVert_{\infty}=\sup\limits_{n}|\mathfrak{m}_{n}|$ on $\ell_{\infty}$, $c$, $c_{0}$,
$\lVert\mathfrak{m}\rVert_{bs}=\sup\bigg|\sum\limits_{n=1}^{j}\mathfrak{m}_{n}\bigg|$ on $cs$ and $bs$ and $\lVert\mathfrak{m}\rVert_{1}=\sum\limits_{n=1}^{\infty}|\mathfrak{m}|$ on $\ell_{1}.$

The forward difference operator $\Delta:\omega\rightarrow\omega$ is defined as follows $\Delta\mathfrak{m}_{n}=\mathfrak{m}_{n}-\mathfrak{m}_{n+1},~n=(1,2,...)$.
 
Hahn \cite{Hahn} was introduced the Hahn space $$h=\{\mathfrak{m}\in\omega:\sum\limits_{n=1}^{\infty}n|\Delta\mathfrak{m}_{n}|<\infty\}\cap c_{0},$$ the Hahn space $h$ is a BK space with $$\lVert\mathfrak{m}\rVert^{'}=\sum\limits_{n=1}^{\infty}n|\Delta\mathfrak{m}_{n}|+\sup\limits_{n}|\mathfrak{m}_{n}|$$ for all $\mathfrak{m}=(\mathfrak{m}_{n})_{n=1}^{\infty}\in h$. 

The Hahn space is a BK-space as well as it qualifies the AK-property with the norm $$\lVert\mathfrak{m}\rVert=\sum\limits_{n=1}^{\infty}n|\Delta\mathfrak{m}_{n}|$$ for all $\mathfrak{m}=(\mathfrak{m}_{n})_{n=1}^{\infty}\in h$, is proved by Rao \cite{Rao}.

According to Goes \cite{Goes}, 
for any complex sequence $d=(d_{n})_{n=1}^{\infty}$ with $d_{n}\neq0$ for all $n$,  the generalized Hahn space $h_{d}$ is defined by 
$$h_{d}=\left\{\mathfrak{m}\in\omega: \sum_{n=1}^{\infty}|d_{n}||\Delta\mathfrak{m}_{n}|<\infty\right\}\cap c_{0},$$
where $\Delta\mathfrak{m}_{n}=\mathfrak{m}_{n}-\mathfrak{m}_{n+1}$ and $\omega$ denotes the set of all complex sequences.

Again the generalized Hahn sequence space was introduced by Malkowsky et al. \cite{Operator} for a sequence $\{d\}_{n=1}^{\infty}$ constitute of positive real numbers, they focused on the case where $d$ is unbounded and monotonically increases sequences.\\
In case, if $d$ is bounded, then it is clear that $$h_{d}=h_{e}=bv_{0}=\{\mathfrak{m}\in\omega: \sum_{n=1}^{\infty}|\mathfrak{m}_{n}-\mathfrak{m}_{n+1}|<\infty\}\cap c_{0}.$$

The space $h_{d}$ is a BK-space and also satisfies the AK-property. Its norm, denoted by $\lVert .\rVert_{h_{d}}$ is defined as follows $$\lVert\mathfrak{m}\rVert_{h_{d}}=\sum\limits_{n=1}^{\infty}d_{n}|\Delta\mathfrak{m}_{n}|<\infty,$$ for all $\mathfrak{m}=(\mathfrak{m}_{n})_{n=1}^{\infty}\in h_{d}$. Various categories of matrix transformation mappings or bounded linear operators from the space $h_{d}$ into the spaces $Y\in\{\ell_{\infty}, c,c_{0},\ell_{1}, h_{d}\}$ as well as transformations from $Y\in\{\ell_{\infty}, c,c_{0},\ell_{1}\}$ into $h_{d}$ are identified and described by the authors. Furthermore, the research explored the norms of the associated bounded linear operators and for those operator classes, they also investigated the Hausdorff measure of noncompactness including an application involving a tridiagonal matrix that defines a Fredholm operator acting on $h_{d}$.

Also in \cite{OTug}, the authors introduced a novel generalized Hahn sequence space, denoted as $h^{p}_{d}$ as follows
$$h^{p}_{d}=\{\mathfrak{m}\in\omega:\sum_{n=1}^{\infty}|d_{n}\Delta\mathfrak{m}_{n}|^{p}<\infty\},$$
where $d=(d_{n})_{n=1}^{\infty}$ represents an unbounded, monotonically increasing sequence of positive real numbers with $d_{n}\neq0$ for every $n\in\mathbb{N}$, and the parameter $p\in(1,\infty)$ and $\Delta\mathfrak{m}_{n}=\mathfrak{m}_{n}-\mathfrak{m}_{n+1}$. Next, they established certain topological properties and demonstrated several inclusion relationships. The $\alpha-$, $\beta-$ and $\gamma-$ duals of $h^{p}_{d}$ are calculated there. In addition, they provided a characterization of the newly defined matrix classes.

\begin{lem}\cite{Mursaleen}
    Let $\mathcal{A}$ be a bounded subset of the space $\mathcal{R}$, where $\mathcal{R}$ is either $\ell_{p}$ $(p\in[1,\infty))$ or $c_{0}$. Consider the projection operator $\mathcal{P}_{k}:\mathcal{R}\rightarrow\mathcal{R}$ defined by $\mathcal{P}_{k}(\mathfrak{m})=(\mathfrak{m_{0}},\mathfrak{m_{1}},...\mathfrak{m_{k}},0,0....)$. Then, the measure of non-compactness $\chi(\mathcal{A})$ can be expressed as:  $$\chi(\mathcal{A})=\lim\limits_{k\rightarrow\infty}\sup\limits_{\mathfrak{m}\in\mathcal{A}}\lVert(\kappa-\mathcal{P}_{k})\mathfrak{m}\rVert,$$
    where the supremum is taken over all $\mathfrak{m}\in\mathcal{A}$ and $I$ denotes the identity operator.
\end{lem}

\section{ Hausdorff MNC on generalized Hahn space}
This part introduces the Hausdorff MNC in the $h_{d}$ space. The following theorem is discussed to reinforce this point.
\begin{thm}
In the Banach space $h_{d}$, consider a bounded subset $\mathcal{A}\subseteq h_{d}$. Then, the Hausdorff measure of non-compactness $\chi$ associated with $\mathcal{A}$ is given by
    \begin{equation}\label{2.1}
     \chi(\mathcal{A}):=\lim\limits_{k\rightarrow\infty}\{\sup\limits_{\mathfrak{m}\in\mathcal{A}}\{\sum\limits_{n\geq k}|d_{n}||\Delta\mathfrak{m}_{n}|\}\}.
    \end{equation}
\end{thm}
\textbf{Proof}: Let us define the operator $\mathcal{P}_{k}:h_{d}\rightarrow h_{d}$ as folllows  $\mathcal{P}_{k}(\mathfrak{m})=(\mathfrak{m}_{1},...\mathfrak{m}_{k},0,0....)$ for $\mathfrak{m}=(\mathfrak{m}_{1}, \mathfrak{m}_{2},...)\in h_{d}$. Then, 
\begin{equation}
    \mathcal{A}\subset\mathcal{P}_{k}\mathcal{A}+(\kappa-\mathcal{P}_{k})\mathcal{A}.
\end{equation}
Then, 
\begin{equation*}
    \begin{split}        \chi({\mathcal{A}})&\leq\chi({\mathcal{P}_{k}\mathcal{A}+(\kappa-\mathcal{P}_{k})\mathcal{A}})\\
        &=\chi((\kappa-\mathcal{P}_{k})\mathcal{A})\\
        &\leq\dim((\kappa-\mathcal{P}_{k})\mathcal{A})\\
        =\sup\limits_{\mathfrak{m}\in\mathcal{A}}\lVert(\kappa-\mathcal{P}_{k})\mathfrak{m}\rVert,
    \end{split}
\end{equation*}
where $\lVert(\kappa-\mathcal{P}_{k})\mathfrak{m}\rVert=\sum\limits_{n=1}^{\infty}|d_{n}||\Delta\mathfrak{m}_{n}|.$ So, when $k$ is sufficiently large. So, 
\begin{equation}\label{2.3}
    \chi({\mathcal{A}})\leq\lim\limits_{k\rightarrow \infty}\sup\limits_{\mathfrak{m}\in\mathcal{L}}\lVert(\kappa-\mathcal{P}_{k})\mathfrak{m}\rVert.
\end{equation}
Conversely, suppose $\varepsilon>0$ and $\{y_{1},y_{2},...y_{j}\}$ is a $[\chi(\mathcal{A})+\varepsilon]$-net of the bounded set $\mathcal{A}$. Then,
$$\mathcal{A}\subset \{y_{1},y_{2},....y_{j}\}+[\chi(\mathcal{A})+\varepsilon]\mathfrak{D}(h_{d}),$$  
for the Banach space $h_{d}$, the closed unit ball is represented as $\mathfrak{D}(h_{d})$. Then, for the bounded set $\mathcal{A}\subset{h_{d}}$, the following inequality holds:
$$\sup\limits_{\mathfrak{m}\in\mathcal{A}}\lVert(\kappa-\mathcal{P}_{k})\mathfrak{m})\rVert\leq\sup\limits_{1\leq i\leq j}\lVert(\kappa-\mathcal{P}_{k})y_{i}\rVert+[\chi(\mathcal{A})+\varepsilon].$$
Consequently, we obtain the limit
\begin{equation}{\label{2.4}}
    \lim\limits_{k\rightarrow \infty}\big(\sup\limits_{\mathfrak{m}\in\mathcal{A}}\lVert(\kappa-\mathcal{P}_{k})\mathfrak{m}\rVert\big)\leq\chi(\mathcal{A})+\varepsilon.
\end{equation}
Since $\varepsilon$ is arbitrary, from (\ref{2.3}) and (\ref{2.4}), the result  (\ref{2.1}) holds. 

\section { Applications}
Next, we show the existence, uniqueness as well as HU-stability for the infinite system (\ref{1.1}) within the generalized Hahn sequence space. To demonstrate the significance of the primary results, two illustrative examples are presented. Examine the following cases.
\begin{enumerate}
    \item [(I)] Let $\varphi_{i}\in\mathfrak{C}(\kappa\times\mathbb{R}^{\infty},~\mathbb{R})$, for each $i\in\mathbb{N}$ be a given family of continuous functions. Define a mapping $\varphi:\kappa\times h_{d}\rightarrow h_{d}$ by 
        $$(\mathfrak{y},~\mathfrak{m})\rightarrow ((\varphi\mathfrak{m})(s))= (\varphi_{1}(s,\mathfrak{m}(s)),~\varphi_{2}(s,\mathfrak{m}(s)),~\varphi_{3}(s,\mathfrak{m}(s)),....),$$ such that the collection of functions $(\varphi(\mathfrak{m})(s))_{s\in\kappa}$ is equicontinuous for $\kappa=[0, 1]$.
    \item[(II)] The following inequalities are satisfied:
    $$|\varphi_{i}(s,~\mathfrak{m}(s))|\leq|\mathfrak{a}_{i}||\mathfrak{m}_{i}|,~~|\varphi_{i}(\mathfrak{n},~\mathfrak{m}(\mathfrak{n}))|\leq|\mathfrak{a}_{i}||\mathfrak{m}_{i}|,$$
    $$|\Delta\varphi_{i}(s,~\mathfrak{m}(s))|\leq|\mathfrak{a}_{i}||\Delta\mathfrak{m}_{i}|,~~|\Delta\varphi_{i}(\mathfrak{n},~\mathfrak{m}(\mathfrak{n}))|\leq|\mathfrak{a}_{i}||\Delta\mathfrak{m}_{i}|.$$ 
    \item [(III)] $\varphi_{i}:\kappa\times\mathbb{R}\rightarrow \mathbb{R}$ is a continuous function such that there exists $\mathfrak{L}>0$ for which 
    $$|\varphi_{i}(s,x(s))-\varphi_{i}(s,\overline{x}(s))|\leq\mathfrak{L}|x-\overline{x}|,$$
    $$|\varphi_{i}(\mathfrak{n},x(\mathfrak{n}))-\varphi_{i}(\mathfrak{n},\overline{x}(\mathfrak{n}))|\leq\mathfrak{L}|x-\overline{x}|.$$
    \item[(IV)] $\forall\xi\in\kappa$ and for any $\mathfrak{m},\overline{\mathfrak{m}}\in\mathbb{R}$, $\exists$ a constant $\mathfrak{L}>0$ such that $$|\varphi(\xi,\mathfrak{m})-\varphi(\xi,{\overline{\mathfrak{m}})}|\leq\mathfrak{L}|\mathfrak{m}-\overline{\mathfrak{m}}|,$$
\end{enumerate}
where $\mathfrak{a}_{i}:\kappa\rightarrow \mathbb{R}$ are continuous and the sequence $(\mathfrak{a}_{i})_{i\in\mathbb{N}}$ is equibounded. 
Put $\mathfrak{A}=\sup\limits_{i\in\mathbb{N}}|\mathfrak{a}_{i}|.$
\begin{thm}\label{3.1}
    Under assumptions (I), (II) and condition $\Bigg(\frac{1}{\Gamma(\beta+1)}-\frac{2}{(2-\mu\varrho^{2})\Gamma(\beta+1)}+\frac{2\mu\varrho^{\beta+1}}{(2-\mu\varrho^{2})\Gamma(\beta+2)}\Bigg)\mathfrak{A}<1$. Equation (\ref{1.1}) has at least one solution $\mathfrak{m}=(\mathfrak{m}_{n})\in\mathfrak{C}(\kappa, h_{d})$ for each $\xi\in\kappa$.
\end{thm}
\textbf{Proof}: Let the operator $\mathfrak{F}: \mathfrak{C}(\kappa,h_{d})\rightarrow \mathfrak{C}(\kappa,h_{d})$ be define as
\begin{equation*}
    \begin{split}
    (\mathfrak{F}\mathfrak{m})(\xi)
    &=\frac{1}{\Gamma(\beta)}\int\limits_{0}^{\xi}(\xi-s)^{\beta-1}\varphi_{i}(s,\mathfrak{m}(s))ds-\frac{2\xi}{(2-\mu\varrho^{2})\Gamma({\beta})}\int\limits_{0}^{1}(1-s)^{\beta-1}\varphi_{i}(s,\mathfrak{m}(s))ds\\
    &+\frac{2\mu\xi}{(2-\mu\varrho^{2})\Gamma({\beta})}\int\limits_{0}^{\varrho}\left( \int\limits_{0}^{s}(s-\mathfrak{n})^{\beta-1}\varphi_{i}(\mathfrak{n},\mathfrak{m}(\mathfrak{n}))d\mathfrak{n}\right)ds.
    \end{split}
\end{equation*}
Additionally, $\mathfrak{C}(\kappa,h_{d})$ is equipped with norm 
$$\lVert \mathfrak{m}\rVert_{\mathfrak{C}(\kappa,h_{d})}=\sum d_{n}|\Delta\mathfrak{m}_{n}|.$$
By applying our assumptions, we obtain
\begin{equation*}
    \begin{split}
    &\lVert(\mathfrak{F}\mathfrak{m})(\xi)\rVert_{h_{d}}\\
     &=\sum_{n=1}^{\infty}d_{n}\bigg|\Delta\bigg(\frac{1}{\Gamma(\beta)}\int\limits_{0}^{\xi}(\xi-s)^{\beta-1}\varphi_{n}(s,\mathfrak{m}(s))d\mathfrak{n}-\frac{2\xi}{(2-\mu\varrho^{2})\Gamma({\beta})}\int\limits_{0}^{1}(1-s)^{\beta-1}\varphi_{n}(s,\mathfrak{m}(s))ds\\
    &+\frac{2\mu\xi}{(2-\mu\varrho^{2})\Gamma({\beta})}\int\limits_{0}^{\varrho}\bigg( \int\limits_{0}^{s}(s-\mathfrak{n})^{\beta-1}\varphi_{n}(\mathfrak{n},\mathfrak{m}(\mathfrak{n}))d\mathfrak{n}\bigg)ds\bigg)\bigg|\\
    &\leq\sum_{n=1}^{\infty}d_{n}\bigg[\frac{1}{\Gamma(\beta)}\int\limits_{0}^{\xi}(\xi-\mathfrak{n})^{\beta-1}|\Delta{\varphi_{n}(s,\mathfrak{m}(s))}ds|-\frac{2\xi}{(2-\mu\varrho^{2})\Gamma({\beta})}\int\limits_{0}^{1}(1-s)^{\beta-1}|\Delta{\varphi_{n}(s,\mathfrak{m}(s))}ds|\\
    &+\frac{2\mu\xi}{(2-\mu\varrho^{2})\Gamma({\alpha})}\int\limits_{0}^{\varrho}\bigg( \int\limits_{0}^{s}(s-\mathfrak{n})^{\beta-1}|\Delta{\varphi_{n}(\mathfrak{n},\mathfrak{m}(\mathfrak{n}))}|d\mathfrak{n}\bigg)ds\bigg]\\
    &\leq\sum_{n=1}^{\infty}d_{n}|\mathfrak{a_{n}}||\Delta \mathfrak{m}_{n}|\bigg(\frac{1}{\Gamma(\beta)}\int\limits_{0}^{\xi}(\xi-n)^{\beta-1}d\mathfrak{n}-\frac{2\xi}{(2-\mu\varrho^{2})\Gamma({\beta})}\int\limits_{0}^{1}(1-s)^{\beta-1}ds\\
    &+\frac{2\mu\xi}{(2-\mu\varrho^{2})\Gamma({\beta})}\int\limits_{0}^{\varrho}\bigg( \int\limits_{0}^{s}(s-\mathfrak{n})^{\beta-1}d\mathfrak{n}\bigg)ds\bigg)\\
    &\leq\Bigg(\frac{\xi^{\beta}}{\Gamma(\beta+1)}-\frac{2\xi}{(2-\mu\varrho^{2})~\Gamma(\beta+1)}+\frac{2\xi\mu\varrho^{\beta+1}}{(2-\mu\varrho^{2})~\Gamma(\beta+2)}\Bigg)\sum_{n=1}^{\infty}d_{n}|\mathfrak{a_{n}}||\Delta \mathfrak{m}_{n}|.
    \end{split}
\end{equation*}
By taking the supremum over the set $\xi$ within the interval $[0,1],$ we get
$$\lVert\mathfrak{F}\mathfrak{m}\rVert_{\mathfrak{C}(\kappa,h_{d})}\leq\Bigg(\frac{1}{\Gamma(\beta+1)}-\frac{2}{(2-\mu\varrho^{2})\Gamma(\beta+1)}+\frac{2\mu\varrho^{\beta+1}}{(2-\mu\varrho^{2})\Gamma(\beta+2)}\Bigg)\mathfrak{A}\lVert \mathfrak{m}\rVert_{\mathfrak{C}(\kappa,h_{d})}.$$
The given inequality can be written as 
\begin{equation}\label{3.1.}
\mathfrak{r}\leq\Bigg(\frac{1}{\Gamma(\beta+1)}-\frac{2}{(2-\mu\varrho^{2})\Gamma(\beta+1)}+\frac{2\mu\varrho^{\beta+1}}{(2-\mu\varrho^{2})\Gamma(\beta+2)}\Bigg)\mathfrak{A}\mathfrak{r},
\end{equation}
where $\mathfrak{r}$ represents the optimal solution of (\ref{3.1.}). Consider
$$\mathfrak{D}=\mathfrak{D}(\xi^{0},\mathfrak{r}_{0})=\{\mathfrak{m}=(\mathfrak{m}_{n})\in\mathfrak{C}(\kappa, h_{d}): \lVert\mathfrak{m}\rVert_{\mathfrak{C}(\kappa,h_{d})}\leq\mathfrak{r}_{0}\}.$$
Clearly, the set $\overline{\mathfrak{D}}=\mathfrak{D}$ is convex, bounded and the function $\varphi$ remains bounded on $\mathfrak{D}$.
Now, suppose $x\in\mathfrak{D}$ and fix $\varepsilon>0.$ Using  the condition (I), there exists $0<\delta$ such that if $\mathfrak{m}\in\mathfrak{D}$ and $\lVert\mathfrak{m}-x\rVert\leq\delta$ then $\lVert(\varphi\mathfrak{m})-(\varphi x)\rVert_{\mathfrak{C}(\kappa,h_{d})}\leq\frac{\varepsilon}{\Bigg(\frac{1}{\Gamma(\beta+1)}-\frac{2}{(2-\mu\varrho^{2})\Gamma(\beta+1)}+\frac{2\mu\varrho^{\beta+1}}{(2-\mu\varrho^{2})\Gamma(\beta+2)}\Bigg)\mathfrak{A}}$.\\
Thus, for every $\xi\in[0,1]$, we have
\begin{equation*}
    \begin{split} 
        &\lVert(\mathfrak{F}\mathfrak{m})(\xi)-(\mathfrak{F}x)(\xi)\rVert_{h_{d}}\\
        &=\sum_{n=1}^{\infty}d_{n}\bigg|\bigg(\frac{1}{\Gamma(\beta)}\int\limits_{0}^{\xi}(\xi-s)^{\beta-1}\Delta\bigg(\varphi_{n}(s,\mathfrak{m}(s))-\varphi_{n}(s,~x(s))\bigg)ds\\
        &-\frac{2\xi}{(2-\mu\varrho^{2})\Gamma({\beta})}\int\limits_{0}^{1}(1-s)^{\beta-1}\Delta\bigg(\varphi_{n}(s,\mathfrak{m}(s))-\varphi_{n}(s,x(s))\bigg)ds\\
        &+\frac{2\mu\xi}{(2-\mu\varrho^{2})\Gamma({\beta})}\int\limits_{0}^{\varrho}\bigg( \int\limits_{0}^{s}(s-\mathfrak{n})^{\beta-1}\Delta\bigg(\varphi_{n}(\mathfrak{n},~\mathfrak{m}(\mathfrak{n}))-\varphi_{n}(\mathfrak{n},~x(\mathfrak{n}))\bigg)d\mathfrak{n}\bigg)ds\bigg|\\
        &\leq\Bigg(\frac{1}{\Gamma(\beta+1)}-\frac{2}{(2-\mu\varrho^{2})\Gamma(\beta+1)}\Bigg)\sum_{n=1}^{\infty}d_{n}|\Delta(\varphi_{n}(s,\mathfrak{m}(s))-\varphi_{n}(s,x(s)))|\\
        &+\Bigg(\frac{2\mu\varrho^{\beta+1}}{(2-\mu\varrho^{2})\Gamma(\beta+2)}\Bigg)\sum_{n=1}^{\infty}d_{n}|\Delta(\varphi_{n}(\mathfrak{n},\mathfrak{m}(\mathfrak{n}))-\varphi_{n}(\mathfrak{n},x(\mathfrak{n})))|\\
        &=\Bigg(\frac{1}{\Gamma(\beta+1)}-\frac{2}{(2-\mu\varrho^{2})\Gamma(\beta+1)}+\frac{2\mu\varrho^{\beta+1}}{(2-\mu\varrho^{2})\Gamma(\beta+2)}\Bigg)\lVert(\varphi\mathfrak{m})-(\varphi x)\rVert_{\mathfrak{C}(\kappa,h_{d})}\\
        &\leq\varepsilon.
    \end{split}
\end{equation*}
Hence $$\lVert(\mathfrak{F}\mathfrak{m})-(\mathfrak{F}x)\rVert_{\mathfrak{C}({\kappa,h_{d})}}\leq\varepsilon.$$
Therefore, the operator $\mathfrak{F}$ is continuous.
We now proceed to show that the function $(\mathfrak{F}\mathfrak{m})$ is continuous on the interval (0,1). Let $\xi_{1}\in(0,1)$ and choose $\xi>\xi_{1}$ such that $|\xi-\xi_{1}|<\varepsilon$, for a given $\varepsilon>0$. 
\begin{align*}
     &\lVert (\mathfrak{F}\mathfrak{m})(\xi)- (\mathfrak{F}\mathfrak{m})(\xi_{1)}\rVert_{h_{d}}\\
     &\leq\sum_{n=1}^{\infty}d_{n}\bigg|\Delta\bigg(\frac{1}{\Gamma(\beta)}\int\limits_{0}^{\xi}(\xi-s)^{\beta-1}\varphi_{n}(s,~\mathfrak{m}(s))ds-\int\limits_{0}^{\xi_{1}}(\xi_{1}-s)^{\beta-1}\varphi_{n}(s,\mathfrak{m}(s))ds\bigg)\bigg|\\
     &-\frac{2(\xi-\xi_{1})}{(2-\mu\varrho^{2})~\Gamma(\beta)}\int\limits_{0}^{1}(1-s)^{\beta-1}\varphi_{n}(s,\mathfrak{m}(s))ds\\
     &+\frac{2\mu\varrho^{\beta+1}(\xi_{1}-\xi)}{(2-\mu\varrho^{2})~\Gamma(\beta)}\int\limits_{0}^{\varrho}\bigg(\int\limits_{0}^{s}(s-\mathfrak{n})^{\beta-1}\varphi_{n}(\mathfrak{n},~\mathfrak{m}(\mathfrak{n}))d\mathfrak{n}\bigg)ds\\
     &\leq\frac{1}{\Gamma(\beta)}\sum_{n=1}^{\infty}d_{n}|\mathfrak{a_{n}}||\Delta \mathfrak{m}_{n}|\bigg[\int\limits_{0}^{\xi}(\xi-s)^{\beta-1}ds-\frac{2(\xi-\xi_{1})}{(2-\mu\varrho^{2})}\int\limits_{0}^{1}(1-s)^{\beta-1}ds\end{align*}\begin{align*}
     &+\frac{2\mu\varrho^{\beta+1}(\xi_{1}-\xi)}{(2-\mu\varrho^{2})}\int\limits_{0}^{\varrho}\bigg( \int\limits_{0}^{s}(s-\mathfrak{n})^{\beta-1}d\mathfrak{n}\bigg)ds\bigg]\\
     &\leq\frac{\mathfrak{A}}{\Gamma(\beta)}\sum_{n=1}^{\infty}d_{n}|\Delta \mathfrak{m}_{n}|\bigg[\bigg(\frac{\xi^{\beta}}{\beta}-\frac{\xi_{1}^{\beta}}{\beta}\bigg)-\frac{2(\xi-\xi_{1})}{(2-\mu\varrho^{2})}\int\limits_{0}^{1}(1-s)^{\beta-1}ds\\
    &+\frac{2\mu\varrho^{\beta+1}(\xi_{1}-\xi)}{(2-\mu\varrho^{2})}\int\limits_{0}^{\varrho}\bigg( \int\limits_{0}^{s}(s-\mathfrak{n})^{\beta-1}d\mathfrak{n}\bigg)ds\bigg].
\end{align*}
Given that $\xi>\xi_{1}$ and $0<\beta\leq1$, it follows that $\frac{\xi^{\beta}-\xi_{1}^{\beta}}{\beta}\leq0$. This confirms the continuity of  $(\mathfrak{F}\mathfrak{m})$ on the interval(0,1).\\
To conclude, we verify that the operator $\mathfrak{F}$ satisfies the conditions outlined in Theorem \ref{1.4}. By applying Proposition \ref{1.5} along with equation (\ref{2.1}), the Hausdorff MNC of a subset $\mathfrak{D}\subset\mathfrak{C}(\kappa, h_{d})$ is expressed as:
$$\chi_{\mathfrak{C}(\kappa,h_{d})}(\mathfrak{D}):=\sup\limits_{\xi\in\kappa}\chi_{h_{d}}(\mathfrak{D}(\xi)),$$
where $\mathfrak{D}(\xi)=\{\mathfrak{m}(\xi); ~\mathfrak{m}\in\mathfrak{D}\}.$
Therefore, we get
\begin{equation*}
    \begin{split}
       & \chi_{h_{d}}(\mathfrak{F}\mathfrak{D})(\xi)\\
        &=\lim\limits_{k\rightarrow\infty}\{\sup\limits_{\mathfrak{m}\in\mathfrak{D}}(\sum\limits_{n\geq k}d_{n}|\Delta(\frac{1}{\Gamma(\beta)}\int\limits_{0}^{\xi}(\xi-s)^{\beta-1}\varphi_{n}(s,\mathfrak{m}(s))ds-\frac{2\xi}{(2-\mu\varrho^{2})\Gamma({\beta})}\int\limits_{0}^{1}(1-s)^{\beta-1}\varphi_{n}(s,\mathfrak{m}(s))ds\\
        &+\frac{2\mu\xi}{(2-\mu\varrho^{2})\Gamma({\beta})}\int\limits_{0}^{\varrho}\bigg( \int\limits_{0}^{s}(s-\mathfrak{n})^{\beta-1}\varphi_{n}(\mathfrak{n},\mathfrak{m}(\mathfrak{n}))d\mathfrak{n}\bigg)ds|)\}\\
        &\leq\lim\limits_{k\rightarrow\infty}\{\sup\limits_{\mathfrak{m}\in\mathfrak{D}}\sum\limits_{n\geq k}d_{n}|\mathfrak{a_{n}}||\Delta\mathfrak{m}_{n}|(\frac{1}{\Gamma(\beta)}\int\limits_{0}^{\xi}(\xi-\mathfrak{n})^{\beta-1}d\mathfrak{n}-\frac{2\xi}{(2-\mu\varrho^{2})\Gamma({\beta})}\int\limits_{0}^{1}(1-\mathfrak{n})^{\beta-1}d\mathfrak{n}\\
        &+\frac{2\mu\xi}{(2-\mu\varrho^{2})\Gamma({\beta})}\int\limits_{0}^{\varrho}\bigg( \int\limits_{0}^{s}(s-\mathfrak{n})^{\beta-1}d\mathfrak{n}\bigg)ds\\
        &\leq\mathfrak{A}\Bigg(\frac{\xi^{\beta}}{\Gamma(\beta+1)}-\frac{2\xi}{(2-\mu\varrho^{2})\Gamma(\beta+1)}+\frac{2\xi\mu\varrho^{\beta+1}}{(2-\mu\varrho^{2})\Gamma(\beta+2)}\Bigg)\lim\limits_{k\rightarrow \infty}\{\sup\limits_{\mathfrak{m}\in\mathfrak{D}}(\sum\limits_{n\geq k}d_{n}|\Delta\mathfrak{m}_{n}|)\}.
    \end{split}
\end{equation*}
Then, we get 
$$\sup\limits_{\xi\in\kappa}\chi_{h_{d}}(\mathfrak{F}\mathfrak{D})(\xi)\leq\mathfrak{A}\Bigg(\frac{1}{\Gamma(\beta+1)}-\frac{2}{(2-\mu\varrho^{2})\Gamma(\beta+1)}+\frac{2\mu\varrho^{\beta+1}}{(2-\mu\varrho^{2})\Gamma(\beta+2)}\Bigg)\chi_{\mathfrak{C}(\kappa, h_{d})}(\mathfrak{D}).$$
This implies that 
\begin{equation}
\chi_{\mathfrak{C}(\kappa, h_{d})}(\mathfrak{F}\mathfrak{D})<\mathfrak{A}\Bigg(\frac{1}{\Gamma(\beta+1)}-\frac{2}{(2-\mu\varrho^{2})\Gamma(\beta+1)}+\frac{2\mu\varrho^{\beta+1}}{(2-\mu\varrho^{2})\Gamma(\beta+2)}\Bigg)\chi_{\mathfrak{C}(\kappa, h_{d})}(\mathfrak{D})<\varepsilon.
\end{equation}
Therefore, 
$$\chi_{\mathfrak{C}(\kappa, h_{d})}(\mathfrak{D})<\frac{\varepsilon}{\mathfrak{A}\Bigg(\frac{1}{\Gamma(\beta+1)}-\frac{2}{(2-\mu\varrho^{2})\Gamma(\beta+1)}+\frac{2\mu\varrho^{\beta+1}}{(2-\mu\varrho^{2})\Gamma(\beta+2)}\Bigg)}.$$
Let
\begin{equation*}
\delta=\varepsilon\left[\frac{1}{{\mathfrak{A}\Bigg(\frac{1}{\Gamma(\beta+1)}-\frac{2}{(2-\mu\varrho^{2})\Gamma(\beta+1)}+\frac{2\mu\varrho^{\beta+1}}{(2-\mu\varrho^{2})\Gamma(\beta+2)}\Bigg)}}-1\right].
\end{equation*}
Therefore, the operator $\mathfrak{F}$ satisfies the Meir-Keeler condensing condition on the set $\mathfrak{D}\subset h_{d}$. By applying Theorem \ref{1.4}, we conclude that $\mathfrak{F}$ admits a fixed point in $\mathfrak{D}$, which guarantees the existence of at least one solution to equation (\ref{1.1}) in the space $\mathfrak{C}(\kappa, h_{d})$.

\section{Uniqueness of the solution}

This part focuses on proving that Eq. (\ref{1.1}) has an unique solution by utilizing the Banach-Point Theorem. The assumptions introduced in Section 4 are applied here to demonstrate the uniqueness.

\begin{thm}\label{4.1}
    Under the assumption (III) with  $\Bigg(\frac{1}{\Gamma(\beta+1)}-\frac{2}{(2-\mu\varrho^{2})\Gamma(\beta+1)}+\frac{2\mu\varrho^{\beta+1}}{(2-\mu\varrho^{2})\Gamma(\beta+2)}\Bigg)\mathfrak{A}<1$. Equation (\ref{1.1}) has a unique solution in $h_{d}$ space.
\end{thm}
\textbf{Proof:}
To proceed, we demonstrate that the operator $\mathcal{H}$ is a contraction. Suppose $x$ and $y$ are solutions of Eq. (\ref{1.1}), then we obtain:
\begin{equation*}
    \begin{split}
&        {\lVert\mathfrak{F}(x)(\xi)-\mathfrak{F}(\overline{x})(\xi)\rVert}_{h_{d}}\\
        &\leq\frac{1}{\Gamma(\beta)}\int\limits_{0}^{\xi}(\xi-s)^{\beta-1}|\varphi_{i}(s,x(s))-\varphi_{i}(s,\overline{x}(s))|ds\\
        &-\frac{2\xi}{(2-\mu\varrho^{2})\Gamma({\beta})}\int\limits_{0}^{1}(1-s)^{\beta-1}|\varphi_{i}(s,x(s))-\varphi_{i}(s,\overline{x}(s))|ds\\
        &+\frac{2\mu\xi}{(2-\mu\varrho^{2})\Gamma({\beta})}\int\limits_{0}^{\varrho}\bigg( \int\limits_{0}^{s}(s-\mathfrak{n})^{\beta-1}|\varphi_{i}(\mathfrak{n},x(\mathfrak{n}))-\varphi_{i}(\mathfrak{n},\overline{x}(\mathfrak{n}))|d\mathfrak{n}\bigg)ds.\\
        &\leq\frac{1}{\Gamma(\beta)}\int\limits_{0}^{\xi}(\xi-s)^{\beta-1}\mathfrak{L}|x-\overline{x}|ds-\frac{2\xi}{(2-\mu\varrho^{2})\Gamma({\beta})}\int\limits_{0}^{1}(1-s)^{\beta-1}\mathfrak{L}|x-\overline{x}|ds\\
        &+\frac{2\mu\xi}{(2-\mu\varrho^{2})\Gamma({\beta})}\int\limits_{0}^{\varrho}\bigg( \int\limits_{0}^{s}(s-\mathfrak{n})^{\beta-1}\mathfrak{L}|x-\overline{x}|d\mathfrak{n}\bigg)ds.\\
        &=\mathfrak{L}\bigg[\frac{1}{\Gamma(\beta)}\int\limits_{0}^{\xi}(\xi-s)^{\beta-1}|x-\overline{x}|ds-\frac{2\xi}{(2-\mu\varrho^{2})\Gamma({\beta})}\int\limits_{0}^{1}(1-s)^{\beta-1}|x-\overline{x}|ds\\
        &+\frac{2\mu\xi}{(2-\mu\varrho^{2})\Gamma({\beta})}\int\limits_{0}^{\varrho}\bigg( \int\limits_{0}^{s}(s-\mathfrak{n})^{\beta-1}|x-\overline{x}|d\mathfrak{n}\bigg)ds\bigg]\\
        &=\mathfrak{L}\bigg[\frac{\xi^{\beta}}{\Gamma(\beta+1)}-\frac{2\xi}{(2-\mu\varrho^{2})\Gamma(\beta+1)}+\frac{2\xi\mu\varrho^{\beta+1}}{(2-\mu\varrho^{2})\Gamma(\beta+2)}\bigg]|x-\overline{x}|.
       \end{split}
\end{equation*}
By taking the $\sup_{\xi\in\kappa}$ on each side of the inequality. Then, we have
\begin{equation*}
    \begin{split}
&        {\lVert\mathfrak{F}(x)(\xi)-\mathfrak{F}(\overline{x})(\xi)\rVert}_{h_{d}}\\
        &\leq\mathfrak{L}\bigg(\frac{1}{\Gamma(\beta+1)}-\frac{2}{(2-\mu\varrho^{2})\Gamma(\beta+1)}+\frac{2\mu\varrho^{\beta+1}}{(2-\mu\varrho^{2})\Gamma(\beta+2)}\Bigg){\lVert x-\overline{x}\rVert}_{h_{d}}.
    \end{split}
\end{equation*}
This implies that 
$$\lVert\mathfrak{F}(x)(\xi)-\mathfrak{F}(\overline{x})(\xi)\rVert_{h_{d}}\leq\mathfrak{L}{\lVert x-\overline{x}\rVert}_{h_{d}}.$$
Thus, the operator $\mathfrak{F}$ satisfies the contraction condition. 
Therefore, we can conclude that $\mathfrak{F}$ possesses a unique fixed point.

\section{ Hyers-Ulam stability analysis}
The stability of a functional equation refers to the existence of an exact solution close to every approximate one. The concept of stability in functional equations was first brought up by Ulam \cite{Ulam} in 1940, asking under what circumstances a linear mapping can be found close to a function is nearly linear. One year later, Hyers \cite{Hyers} resolved this problem for approximately additive mappings in Banach spaces. In 1978, Rassias \cite{Rassias} extended Hyers result by demonstrating that unique linear mappings exist close to approximately additive mappings. Over time, numerous researcher expanded on these findings in different directions. Ulam stability is a term that encompasses both Hyers-Ulam-Rassias stability. This concept later prompted researchers to explore the Ulam stability of differential equations of integer order.
Wang et al. \cite{wang} conducted a study on the Ulam stability of a first-order differential equation under a given boundary value condition. Even though boundary value problems involving fractional integrals are, they are difficult to solve and exact solutions are rarely attainable. Therefore, analyzing the Ulam stability of fractional differential equations is essential, as it enables the substitution of exact solutions with approximate ones. Researchers have thus far investigated the Ulam stability, as well as the existence and uniqueness of solutions, for fractional differential equations involving Hilfer-Hadamard, Caputo and Caputo-Fabrizio fractional derivatives.

This section explores the concept of HU-stability for the proposed  Equation \eqref{1.1}.
\begin{defn}
    The nonlinear fractional differential equation \eqref{1.1} satisfies the UHS if $\exists$ a number $\mathcal{G}$ s.t. for every $\varepsilon>0$ 
    \begin{equation}\label{5.1}
      \left|^c\mathfrak{D}^{\beta}\mathfrak{m}(\xi)-\varphi(\xi,~\mathfrak{m}(\xi)) \right|\leq\varepsilon,
    \end{equation} 
    $\exists$ a solution $y\in h_{d}$ of \eqref{1.1} in such a way that
$$\left|\mathfrak{m}(\xi)-y(\xi)\right|\leq\mathcal{G}\varepsilon,~\xi\in\kappa.$$
Such $\mathcal{G}$ is termed as a HU-stability constant.
\end{defn}
\begin{rem}\label{1.11}
    A function $\mathfrak{m}\in h_{d}$ is a solution for the inequality \eqref{5.1} if and only if there exists a function $\mathfrak{g}\in h_{d}$ (depending on $\mathfrak{m}$) such that 
    \item{i.} $|\mathfrak{g}(\xi)|\leq\varepsilon, ~\forall\xi\in\kappa,$ 
    \item{ii.} $^c\mathfrak{D}^{\beta}\mathfrak{m}(\xi)=\varphi(\xi,~\mathfrak{m}(\xi))+\mathfrak{g}(\xi),~\xi\in\kappa.$
\end{rem}
\begin{lem}
    Let $\beta\in(0,1]$. If $\mathfrak{m}\in h_{d}$ is a solution of the equation \eqref{5.1}, then it must also fulfill the following inequality 
    $$\left|\mathfrak{m}(\xi)-y(\xi)\right|\leq\mathcal{G}\varepsilon.$$
\end{lem}
\textbf{Proof}: Let $\mathfrak{m}$ denotes the solution of the inequality \eqref{5.1}, considering the Remark \ref{1.11}, we can conclude that 
\begin{equation}\label{5.2}
    \begin{split}
        ^c\mathfrak{D}^{\beta}\mathfrak{m}(\xi)=\varphi(\xi,~\mathfrak{m}(\xi))+\mathfrak{g}(\xi),\\
        \mathfrak{m}(0)=0,~\mathfrak{m}(1)=\mu\int\limits_{0}^{\varrho}\mathfrak{m}(s)ds,~ 0<\varrho<1.
    \end{split}
\end{equation}
Therefore, the solution of \eqref{5.2} becomes of the following form 
\begin{equation*}
        \begin{split}
           \mathfrak{m}(\xi)
           &=\frac{1}{\Gamma(\beta)}\int\limits_{0}^{\xi}(\xi-s)^{\beta-1}\varphi(s,\mathfrak{m}(s))ds-\frac{2\xi}{(2-\mu\varrho^{2})\Gamma({\beta})}\int\limits_{0}^{1}(1-s)^{\beta-1}\varphi(s,\mathfrak{m}(s))ds\\
           &+\frac{2\mu\xi}{(2-\mu\varrho^{2})\Gamma({\beta})}\int\limits_{0}^{\varrho}\left( \int\limits_{0}^{s}(s-\mathfrak{n})^{\beta-1}\varphi(\mathfrak{n},\mathfrak{m}(\mathfrak{n})d\mathfrak{n}\right)ds\\
           &+\frac{1}{\Gamma(\beta)}\int\limits_{0}^{\xi}(\xi-s)^{\beta-1}\mathfrak{g}(s)ds-\frac{2\xi}{(2-\mu\varrho^{2})\Gamma({\beta})}\int\limits_{0}^{1}(1-s)^{\beta-1}\mathfrak{g}(s)ds\\
           &+\frac{2\mu\xi}{(2-\mu\varrho^{2})\Gamma({\beta})}\int\limits_{0}^{\varrho}\left( \int\limits_{0}^{s}(s-\mathfrak{n})^{\beta-1}\mathfrak{g}(s)d\mathfrak{n}\right)ds.
        \end{split}
    \end{equation*}
To simplify, let us define denote $y(\xi)$ to represent the sum of the terms excluding $\mathfrak{g}$, this allows us to rewrite the expression as follows: 
\begin{equation*}
    \begin{split}
        y(\xi)
        &=\frac{1}{\Gamma(\beta)}\int\limits_{0}^{\xi}(\xi-s)^{\beta-1}\varphi(s,\mathfrak{m}(s))ds-\frac{2\xi}{(2-\mu\varrho^{2})\Gamma({\beta})}\int\limits_{0}^{1}(1-s)^{\beta-1}\varphi(s,\mathfrak{m}(s))ds\\
           &+\frac{2\mu\xi}{(2-\mu\varrho^{2})\Gamma({\beta})}\int\limits_{0}^{\varrho}\left( \int\limits_{0}^{s}(s-\mathfrak{n})^{\beta-1}\varphi(\mathfrak{n},\mathfrak{m}(\mathfrak{n})d\mathfrak{n}\right)ds.
    \end{split}
\end{equation*}
So, we obtain 
\begin{equation*}
    \begin{split}
        |\mathfrak{m}(\xi)-y(\xi)|
        &\leq\frac{1}{\Gamma(\beta)}\int\limits_{0}^{\xi}(\xi-s)^{\beta-1}|\mathfrak{g}(s)|ds-\frac{2\xi}{(2-\mu\varrho^{2})\Gamma({\beta})}\int\limits_{0}^{1}(1-s)^{\beta-1}|\mathfrak{g}(s)|ds\\
        &+\frac{2\mu\xi}{(2-\mu\varrho^{2})\Gamma({\beta})}\int\limits_{0}^{\varrho}\left( \int\limits_{0}^{s}(s-\mathfrak{n})^{\beta-1}|\mathfrak{g}(s)|d\mathfrak{n}\right)ds\\
        &\leq\bigg[\frac{\xi^{\beta}}{\Gamma(\beta+1)}-\frac{2\xi}{(2-\mu\varrho^{2})\Gamma(\beta+1)}+\frac{2\xi\mu\varrho^{\beta+1}}{(2-\mu\varrho^{2})\Gamma(\beta+2)}\bigg]\varepsilon,
    \end{split}
\end{equation*}
i.e., $$|\mathfrak{m}(\xi)-y(\xi)|\leq\mathcal{G}\varepsilon.$$
\begin{thm}
    Under the hypothesis (III) and (IV) with $\mathcal{G}\mathfrak{L}<1$, the equation \eqref{1.1} is UHS.
\end{thm}
\textbf{Proof:} Let us suppose $z$ be the unique solution of the equation 
\begin{equation}\label{5.3}
    \begin{split}
        ^c\mathfrak{D}^{\beta}y(\xi)=\varphi(\xi,~z(\xi)),~0<\xi<1,\\
        with ~B.C.~
        z(0)=0,~z(1)=\mu\int\limits_{0}^{\varrho}z(s)ds,~ 0<\varrho<1.
    \end{split}
\end{equation}
Therefore, the solution of (\ref{5.3}) for $\mathfrak{m}\in\kappa$ is 
\begin{equation*}
        \begin{split}
           \mathfrak{m}(\xi)
           &=\frac{1}{\Gamma(\beta)}\int\limits_{0}^{\xi}(\xi-s)^{\beta-1}\varphi(s,\mathfrak{m}(s))ds-\frac{2\xi}{(2-\mu\varrho^{2})\Gamma({\beta})}\int\limits_{0}^{1}(1-s)^{\beta-1}\varphi(s,\mathfrak{m}(s))ds\\
           &+\frac{2\mu\xi}{(2-\mu\varrho^{2})\Gamma({\beta})}\int\limits_{0}^{\varrho}\left( \int\limits_{0}^{s}(s-\mathfrak{n})^{\beta-1}\varphi(\mathfrak{n},\mathfrak{m}(\mathfrak{n})d\mathfrak{n}\right)ds.
        \end{split}
    \end{equation*}
Suppose
\begin{equation}\label{5.4}
|\mathfrak{m}(\xi)-z(\xi)|\leq|\mathfrak{m}(\xi)-y(\xi)|+|y(\xi)-z(\xi)|.
\end{equation}
So, by using Lemma \ref{5.2} in (\ref{5.4}), we have
\begin{equation}
    \begin{split}
        &|\mathfrak{m}(\xi)-z(\xi)|\leq\mathcal {G}\varepsilon
        +\frac{1}{\Gamma(\beta)}\int\limits_{0}^{\xi}(\xi-s)^{\beta-1}|\varphi(s, y(s))ds-\varphi(s,z(s))|\\
        &-\frac{2\xi}{(2-\mu\varrho^{2})\Gamma({\beta})}\int\limits_{0}^{1}(1-s)^{\beta-1}|\varphi(s,y(s))-\varphi(s,z(s))|ds\\
        &+\frac{2\mu\xi}{(2-\mu\varrho^{2})\Gamma({\beta})}\int\limits_{0}^{\varrho}\left( \int\limits_{0}^{s}(s-\mathfrak{n})^{\beta-1}|\varphi(\mathfrak{n},y(\mathfrak{n})d\mathfrak{n}-\varphi(\mathfrak{n},z(\mathfrak{n})|\right)ds.
    \end{split}
\end{equation}
From the assumption $[IV]$, we have  
$$|\varphi(\xi,\mathfrak{m}(\xi))-\varphi(\xi,{\overline{\mathfrak{m}}(\xi))}|\leq\mathfrak{L}|\mathfrak{m}(\xi)-\overline{\mathfrak{m}}(\xi)|.$$
So, we have
$$\lVert \mathfrak{m}-z\rVert_{h_{d}}\leq\frac{\mathcal{G}\varepsilon}{1-\mathcal{G}\mathfrak{L}}=\mathcal{G}_{0}\varepsilon,$$
where $$\mathcal{G}_{0}=\frac{\mathcal{G}}{1-\mathcal{G}\mathfrak{L}}.$$
So that $$\mathcal{G}\mathfrak{L}<1.$$

\section{ Examples}
To verify the accuracy of our results, we present multiple examples in this section.
\begin{exm} Let us examine the following equation 
\begin{equation}\label{2}
    \begin{split}
        ^c\mathfrak{D}^{\frac{1}{4}}\mathfrak{m}_{i}(\xi)=\sum\limits_{i=1}^{\infty}\frac{e^{-4s}}{(i+3)^{2}}\cos(s+e^{2s})\mathfrak{m}_{i}(s)\\
        \xi\in[0,1],~\beta=\frac{1}{4},~\varrho=\frac{1}{3},~ \mu=1.\\
        \mathfrak{m}(0)=0,~\mathfrak{m}(1)=\int\limits_{0}^{\frac{1}{3}}\mathfrak{m}(s)ds,~ 0<\varrho<1.
    \end{split}
\end{equation}
Comparing Eq. (\ref{2}) with Eq. (\ref{1.1}), we have
$$\varphi_{i}(\xi,~\mathfrak{m}(\xi))=\sum\limits_{i=1}^{\infty}\frac{e^{-4s}}{(i+3)^{2}}\cos(s+e^{2s})\mathfrak{m}_{i}(s).$$
Clearly, $\sum\limits_{i=1}^{\infty}\frac{e^{-4s}}{(i+3)^{2}}\cos(s+e^{2s})\mathfrak{m}_{i}(s)$ is continuous on $\kappa\in[0,1].$
Note that, for any $\xi\in\kappa$, if $\mathfrak{m}(\xi)=(\mathfrak{m}_{i}(\xi))\in h_{d},$ then $(\varphi_{i}(\xi,~\mathfrak{m}(\xi)))\in h_{d}.$ So, by taking $x(\xi)=(x_{i}(\xi))\in h_{d}$ with $\Vert\mathfrak{m}(\xi)-x(\xi)\rVert_{h_{d}}\leq \delta(\varepsilon)=9\varepsilon,$ we have
$$\lVert\varphi_{i}(\xi,~\mathfrak{m}(\xi))-\varphi_{i}(\xi,~x(\xi))\rVert_{h_{d}}\leq\frac{1}{9}\lVert\mathfrak{m}(\xi)-x(\xi)\rVert_{h_{d}}=\varepsilon,$$ 
This confirms that condition (I) is satisfied. Next, we consider condition (II), for which we have:
\begin{equation*}
\begin{split}
|\varphi_{i}(\xi,~\mathfrak{m}(\xi))|=&|\sum\limits_{i=1}^{\infty}\frac{e^{-4s}}{(i+3)^{2}}\cos(s+e^{2s})\mathfrak{m}_{i}(s)|\\
&\leq|\mathfrak{a}_{i}||\mathfrak{m}_{i}|.
\end{split}
\end{equation*}
Also
\begin{equation*}
\begin{split}
|\Delta\varphi_{i}(\xi,~\mathfrak{m}(\xi))|=&|\Delta\sum\limits_{i=1}^{\infty}\frac{e^{-4s}}{(i+3)^{2}}\cos(s+e^{2s})\mathfrak{m}_{i}(s)|\\
&\leq|\mathfrak{a}_{i}||\Delta\mathfrak{m}_{i}|.
\end{split}
\end{equation*}
$\mathfrak{a}_{i}(s)=\frac{e^{-4s}}{9}$ are continuous and $(\mathfrak{a}_{i}(\xi))_{i\in\mathbb{N}}$ is equibounded, by $\mathfrak{A}\leq\frac{1}{9}$ and 
$$\Bigg(\frac{1}{\Gamma(\beta+1)}-\frac{2}{(2-\mu\varrho^{2})\Gamma(\beta+1)}+\frac{2\mu\varrho^{\beta+1}}{(2-\mu\varrho^{2})\Gamma(\beta+2)}\Bigg)\mathfrak{A}=0.370<1.$$
Then, Theorem \ref{3.1} grantees that equation (\ref{2}) possesses at least one solution in the space $\mathfrak{C}([0,1], h_{d}).$
\end{exm}
\begin{exm} We consider 
    \begin{equation}\label{4.2}
    \begin{split}
        ^c\mathfrak{D}^{\frac{1}{5}}\mathfrak{m}_{i}(\xi)
        &=\sum\limits_{l=1}^{\infty}\arctan\bigg(\frac{5}{1+(5l+2)(5l-3)}\bigg)e^{-3s}\ln\big(\frac{8}{9}+|m_{i}(\xi)|\big)\\
        &\mathfrak{m}(0)=0,~\mathfrak{m}(1)=\frac{1}{2}\int\limits_{0}^{\frac{1}{6}}\mathfrak{m}(s)ds,~ 0<\varrho<1.\\
        &\xi\in[0,1],~\beta=\frac{1}{5},~\varrho=\frac{1}{6},~ \mu=\frac{1}{2}.
    \end{split}
    \end{equation}
    Here
$$\varphi_{i}(\xi,~\mathfrak{m}(\xi))=\sum\limits_{l=1}^{\infty}\arctan\bigg(\frac{5}{1+(5l+2)(5l-3)}\bigg)e^{-3s}\ln\big(\frac{8}{9}+|m_{i}(\xi)|\big).$$
It is clear that $\sum\limits_{l=1}^{\infty}\arctan\bigg(\frac{5e^{-3s}}{1+(5l+2)(5l-3)}\bigg)\ln\big(\frac{8}{9}+|m_{i}(\xi)|\big)$ is continuous on $\kappa\in[0,1]$. Keep in mind that, for any $\xi\in\kappa$, if $\mathfrak{m}(\xi)=(\mathfrak{m}_{i}(\xi))\in h_{d},$ then $(\varphi_{i}(\xi,~\mathfrak{m}(\xi)))\in h_{d}.$ Therefore, by taking $x(\xi)=(x_{i}(\xi))\in h_{d}$ with $\Vert\mathfrak{m}(\xi)-x(\xi)\rVert_{h_{d}}\leq \delta(\varepsilon)=3\varepsilon,$ we have
$$\lVert\varphi_{i}(\xi,~\mathfrak{m}(\xi))-\varphi_{i}(\xi,~x(\xi))\rVert_{h_{d}}\leq\frac{1}{3}\lVert\mathfrak{m}(\xi)-x(\xi)\rVert_{h_{d}}=\varepsilon.$$ The condition (I) is fulfilled.
Now 
\begin{equation*}
\begin{split}
|\varphi_{i}(\xi,~\mathfrak{m}(\xi))|=&|\sum\limits_{l=1}^{\infty}\arctan\bigg(\frac{5}{1+(5l+2)(5l-3)}\bigg)e^{-3s}\ln\big(\frac{8}{9}+|m_{i}(\xi)|\big)|\\
&\leq|\mathfrak{a}_{i}||\mathfrak{m}_{i}|.
\end{split}
\end{equation*}
Also 
\begin{equation*}
\begin{split}
|\Delta\varphi_{i}(\xi,~\mathfrak{m}(\xi))|=&|\Delta\sum\limits_{l=1}^{\infty}\arctan\bigg(\frac{5}{1+(5l+2)(5l-3)}\bigg)e^{-3s}\ln\big(\frac{8}{9}+|m_{i}(\xi)|\big)|\\
&\leq|\mathfrak{a}_{i}||\Delta\mathfrak{m}_{i}|.
\end{split}
\end{equation*}
Here, $\mathfrak{a}_{i}(s)=\frac{e^{-3s}}{3}$ are continuous and $(\mathfrak{a}_{i}(\xi))_{i\in\mathbb{N}}$ is equibounded with $\mathfrak{A}\leq\frac{1}{3}$. So, condition (II) is also satisfied.
$$\Bigg(\frac{1}{\Gamma(\beta+1)}-\frac{2}{(2-\mu\varrho^{2})\Gamma(\beta+1)}+\frac{2\mu\varrho^{\beta+1}}{(2-\mu\varrho^{2})\Gamma(\beta+2)}\Bigg)\mathfrak{A}=0.016<1.$$
As a result, by applying Theorem \ref{3.1}, we conclude that equation (\ref{2}) possesses at least one solution in the space $\mathfrak{C}([0,1], h_{d}).$
\end{exm}
\section*{Declaration}
\begin{itemize}
	\item Funding$\colon$ This study was conducted without receiving any funding
	\item Author's Contribution$\colon$ All author's contributed equally.
	\item Conflict of Interest$\colon$ The author declares that they have no conflict of interest.
	\item Data Availability$\colon$ Not applicable.
\end{itemize}
\section*{Acknowledgment} 
 The research of the first author (Khurshida Parvin) is supported by SERB, DST, New Delhi, India, under the application reference number DST/INSPIRE/03/2022/001487.

\end{document}